\documentclass{amsart}
\title  {On the linear independency of monoidal natural transformations}
\author {Kenichi Shimizu}
\date   {}

\usepackage{amsmath,amssymb,amsthm,amscd}
\theoremstyle{plain}
\newtheorem{lemma}              {Lemma}[section]
\newtheorem{proposition}        [lemma]{Proposition}
\newtheorem{theorem}            [lemma]{Theorem}
\newtheorem{corollary}          [lemma]{Corollary}

\theoremstyle{definition}

\theoremstyle{remark}
\newtheorem{remark}             [lemma]{Remark}

\numberwithin{equation}{section}

\newcommand{\id}        {\mathop{\mathrm{id}}\nolimits}
\newcommand{\bfk}       {\mathbf{k}}
\newcommand{\Hom}       {\mathop{\mathrm{Hom}}\nolimits}
\newcommand{\End}       {\mathop{\mathrm{End}}\nolimits}
\newcommand{\Nat}       {\mathop{\mathrm{Nat}}\nolimits}
\newcommand{\Ker}       {\mathop{\mathrm{Ker}}\nolimits}
\newcommand{\Aut}       {\mathop{\mathrm{Aut}}\nolimits}
\newcommand{\Rep}       {\mathop{\mathrm{Rep}}\nolimits}
\newcommand{\Piv}       {\mathop{\mathrm{Piv}}\nolimits}
\newcommand{\Map}       {\mathop{\mathrm{Map}}\nolimits}

\newenvironment{enumalph}
{\begin{enumerate}\renewcommand{\labelenumi}{\textnormal{\hbox to 1.4em{(\hfil\alph{enumi}\hfil)}}}}
  {\end{enumerate}}

\begin{document}

\maketitle

\begin{abstract}
  Let $F, G: \mathcal{I} \to \mathcal{C}$ be strong monoidal functors from a skeletally small monoidal category $\mathcal{I}$ to a tensor category $\mathcal{C}$ over an algebraically closed field $\bfk$. The set $\Nat(F, G)$ of natural transformations $F \to G$ is naturally a vector space over $\bfk$. We show that the set $\Nat_\otimes(F, G)$ of monoidal natural transformations $F \to G$ is linearly independent as a subset of $\Nat(F, G)$.

  As a corollary, we can show that the group of monoidal natural automorphisms on the identity functor on a finite tensor category is finite. We can also show that the set of pivotal structures on a finite tensor category is finite.
\end{abstract}

\section{Introduction}

Monoidal categories \cite{MR1712872} arise in many contexts in mathematics. In this note, we prove a basic fact on monoidal natural transformations between strong monoidal functors. To describe our result, we recall some definitions. Throughout, we work over an algebraically closed field $\bfk$. By a {\em tensor category} over $\bfk$ we mean a $\bfk$-linear abelian rigid monoidal category $\mathcal{C}$ satisfying the following conditions:
\begin{itemize}
\item The unit object $\mathbf{1} \in \mathcal{C}$ is simple.
\item The tensor product $\otimes: \mathcal{C} \times \mathcal{C} \to \mathcal{C}$ is $\bfk$-linear in both variables.
\item Every object of $\mathcal{C}$ is of finite length.
\item Every hom-set in $\mathcal{C}$ is finite-dimensional over $\bfk$.
\end{itemize}

Let $\mathcal{C}$ be a tensor category over $\bfk$ and let $\mathcal{I}$ be a skeletally small monoidal category. If $F, G: \mathcal{I} \to \mathcal{C}$ are functors, the set $\Nat(F, G)$ of natural transformations $F \to G$ is naturally a vector space over $\bfk$. Now we suppose that both $F$ and $G$ are strong monoidal functors. Then we can consider the set $\Nat_\otimes(F, G)$ of monoidal natural transformations. Our main result in this note is the following:

\begin{theorem}
  \label{thm:lin-indep}
  $\Nat_\otimes(F, G) \subset \Nat(F, G)$ is linearly independent.
\end{theorem}

A $\bfk$-linear abelian category $\mathcal{A}$ is said to be {\em finite} if it is $k$-linearly equivalent to the category of fi\-nite\-ly-gen\-er\-at\-ed modules over a fi\-nite-di\-men\-sion\-al $\bfk$-algebra. We present some applications of Theorem~\ref{thm:lin-indep} to finite tensor categories \cite{MR2119143}.

By a {\em tensor functor} we mean a $k$-linear strong monoidal functor between tensor categories. Let $F, G: \mathcal{C} \to \mathcal{D}$ be right exact tensor functors from a finite tensor category $\mathcal{C}$ to a tensor category $\mathcal{D}$. Then it turns out that $\Nat(F, G)$ is fi\-nite-di\-men\-sion\-al (see Lemma~\ref{lem:nat-fin-dim}). Hence:

\begin{corollary}
  \label{cor:fin-nat}
  Under the above assumptions, $\Nat_\otimes(F, G)$ is finite.
\end{corollary}

Let $\Aut_\otimes(F) \subset \Nat_\otimes(F, F)$ denote the group of monoidal natural automorphism on a monoidal functor $F: \mathcal{C} \to \mathcal{D}$. As an immediate consequence of Corollary~\ref{cor:fin-nat}, we have the following:

\begin{corollary}
  \label{cor:fin-auto}
  Let $\mathcal{C}$ be a finite tensor category and let $F: \mathcal{C} \to \mathcal{C}$ be a right exact tensor functor. Then $\Aut_\otimes(F)$ is finite. In particular, $\Aut_\otimes(\id_\mathcal{C})$ is finite.
\end{corollary}

We give some remarks on the structure of $\Aut_\otimes(\id_\mathcal{C})$ in Section~\ref{sec:structure-auto}. We note that the finiteness of $\Aut_\otimes(\id_\mathcal{C})$ is well-known in the case where $\mathcal{C}$ is the category $\Rep(H)$ of fi\-nite-di\-men\-sion\-al representations of a fi\-nite-di\-men\-sion\-al Hopf algebra $H$. In fact, then $\Aut_\otimes(\id_\mathcal{C})$ is isomorphic to the group of central group\-like elements of $H$.

There is another important corollary of Corollary~\ref{cor:fin-nat}. A {\em pivotal structure} on a rigid monoidal category $\mathcal{C}$ is an element of $\Piv(\mathcal{C}) := \Nat_\otimes(\id_\mathcal{C}, (-)^{* *})$, where $(-)^*: \mathcal{C} \to \mathcal{C}$ is the left duality functor.

\begin{corollary}
  \label{cor:fin-piv}
  If $\mathcal{C}$ is a finite tensor category, then $\Piv(\mathcal{C})$ is finite.
\end{corollary}

Also this corollary is well-known in the case where $\mathcal{C} = \Rep(H)$ for some fi\-nite-di\-men\-sion\-al Hopf algebra $H$. In fact, then $\Piv(\mathcal{C})$ is in one-to-one correspondence between group\-like elements $g \in H$ such that $S^2(x) = g x g^{-1}$ for all $x \in H$, where $S$ is the antipode of $H$.

\begin{remark}
  We should remark that every finite tensor category is equivalent as a tensor category to the category of fi\-nite-di\-men\-sion\-al representations of a weak quasi-Hopf algebra \cite{MR2119143}. It will be interesting to interpret our results in terms of ``group\-like elements'' of weak quasi-Hopf algebras ({\it cf.} \cite{MR1932332}).  
\end{remark}

\section{Proof of the main theorem and its corollary}

\subsection{Proof of Theorem~\ref{thm:lin-indep}}

Let $\mathcal{C}$ be a tensor category over $\bfk$. Without loss of generality, we may suppose $\mathcal{C}$ to be strict. Let $V \in \mathcal{C}$ be an object. The tensor product of $\mathcal{C}$ defines two $\bfk$-linear endo\-functors $V \otimes(-)$ and $(-) \otimes V$ on $\mathcal{C}$. Both these functors are exact since $\mathcal{C}$ is rigid \cite{MR1797619}.

\begin{lemma}
  \label{lem:faithful-1}
  If $V \ne 0$, then $V \otimes (-)$ and $(-) \otimes V$ are faithful.
\end{lemma}
\begin{proof}
  Suppose that $V \ne 0$. Then the evaluation morphism $d_V: V^* \otimes V \to \mathbf{1}$ is an epimorphism. Indeed, otherwise, since $\mathbf{1}$ is simple, $d_V = 0$. It follows from the rigidity axiom that $\id_V = 0$, and hence $V = 0$, a contradiction.

  Let $X, Y \in \mathcal{C}$ be objects. Consider the linear map
  \begin{equation*}
    (-)^\natural: \Hom_\mathcal{C}(V \otimes X, V \otimes Y) \to \Hom_\mathcal{C}(V^* \otimes V \otimes X, Y)
  \end{equation*}
  given by $\varphi^\natural = (d_V \otimes \id_Y) \circ (\id_{V^*} \otimes \varphi)$. If $f: X \to Y$ is a morphism in $\mathcal{C}$, then
  \begin{equation*}
    (\id_V \otimes f)^\natural = f \circ (d_V \otimes \id_X).
  \end{equation*}

  If $\id_V \otimes f = 0$, then, by the above equation, we have $f \circ (d_V \otimes \id_X) = 0$. As we observed above, $d_V$ is an epimorphism. Since the tensor product of $\mathcal{C}$ is exact, also $d_V \otimes \id_X$ is an epimorphism. Therefore, we conclude that $f = 0$. This means that the functor $V \otimes (-)$ is faithful.

  The faithfulness of $(-) \otimes V$ can be proved in a similar way.
\end{proof}

Let $\mathcal{I}$ be a skeletally small category and let $F, G: \mathcal{I} \to \mathcal{C}$ be functors. Then the set $\Nat(F, G)$ of natural transformations $F \to G$ is a vector space over $\bfk$. Let $V, W \in \mathcal{C}$ be objects. The tensor product of $\mathcal{C}$ defines a linear map
\begin{equation*}
  \Hom_\mathcal{C}(V, W) \otimes_\bfk \Nat(F, G) \to \Nat(V \otimes F(-), W \otimes G(-)).
\end{equation*}

\begin{lemma}
  \label{lem:faithful-2}
  The above map is injective.
\end{lemma}

\begin{proof}
  Let $f_1, \cdots, f_m \in \Nat(F, G)$ be linearly independent elements. It suffices to show that if $c_1, \cdots, c_m: V \to W$ are morphisms in $\mathcal{C}$ such that
  \begin{equation}
    \label{eq:lin-indep-1}
    c_1 \otimes f_1|_X + \cdots + c_m \otimes f_m|_X = 0
  \end{equation}
  for all $X \in \mathcal{I}$, then $c_i = 0$ for $i = 1, \cdots, m$.

  We show the above claim by induction on the length $\ell(V)$ of $V$. If $\ell(V) = 0$, our claim is obvious. Suppose that $\ell(V) \ge 1$. Then there exists a simple subobject of $V$, say $L$. Let $K$ be the image of the morphism
  \begin{equation*}
    (c_1|_L, c_2|_L, \cdots, c_m|_L): L^{\oplus m} \to W.
  \end{equation*}
  Since $L$ is simple, $K \cong L^{\oplus n}$ for some $n \ge 0$. Let $p_j: K \cong L^{\oplus n} \to L$ be the $j$-th projection. Then $p_j \circ c_i|_L = \lambda_{i j} \id_L$ for some $\lambda_{i j} \in \bfk$. By~\eqref{eq:lin-indep-1},
  \begin{equation*}
    \id_L \otimes \left( \sum_{i = 1}^m \lambda_{i j} f_i|_X \right)
    = (p_j \otimes \id_{G(X)}) \circ \left( \sum_{i = 1}^m c_i|_L \otimes f_i|_X \right) = 0.
  \end{equation*}
  Lemma~\ref{lem:faithful-1} yields that $\sum_{i = 1}^m \lambda_{i j} f_i|_X = 0$ for all $X \in \mathcal{I}$. By the linear independence of $f_i$'s, we have $\lambda_{i j} = 0$ for all $i$ and $j$. This means that $c_i|_L = 0$ for all $i$.

  Let $p: V \to V/L$ be the projection. By the above observation, for each $i$, there exists a morphism $\overline{c}_i: V/L \to W$ such that $c_i = \overline{c}_i \circ p$. By~\eqref{eq:lin-indep-1},
  \begin{equation*}
    \sum_{i = 1}^m \overline{c}_i \otimes f_i|_X
    = \left( \sum_{i = 1}^m c_i \otimes f_i|_X\right) \circ (p \otimes \id_{F(X)}) = 0
  \end{equation*}
  for all $X \in \mathcal{I}$. Note that $\ell(V/L) = \ell(V) - 1$. By the induction hypothesis, $\overline{c}_i = 0$ for all $i$. Therefore, $c_i = \overline{c}_i \circ p = 0$.
\end{proof}

\begin{remark}
  In the case where $F$ and $G$ are constant functors sending all objects of $\mathcal{I}$ respectively to $X \in \mathcal{C}$ and $Y \in \mathcal{C}$, Lemma~\ref{lem:faithful-2} states that the map
  \begin{equation*}
    \Hom_\mathcal{C}(V, W) \otimes_\mathbf{k} \Hom_\mathcal{C}(X, Y)
    \to \Hom_\mathcal{C}(V \otimes X, W \otimes Y)
  \end{equation*}
  induced from the tensor product is injective.
\end{remark}

Now we generalize Lemma~\ref{lem:faithful-2}. Let $\mathcal{I}_i$ ($i = 1, 2$) be skeletally small categories and let $F_i: \mathcal{I}_i \to \mathcal{C}$ ($i = 1, 2$) be functors. We denote by $F_1 \otimes F_2$ the functor
\begin{equation*}
  \begin{CD}
    F_1 \otimes F_2: \mathcal{I}_1 \times \mathcal{I}_2
    @>{F_1 \times F_2}>> \mathcal{C} \times \mathcal{C}
    @>{\otimes}>> \mathcal{C}.
  \end{CD}
\end{equation*}

\begin{lemma}
  \label{lem:faithful-3}
  Let $F_i, G_i: \mathcal{I}_i \to \mathcal{C}$ $(i = 1, 2)$ be functors. The map
  \begin{equation*}
    \Nat(F_1, G_1) \otimes_\bfk \Nat(F_2, G_2) \to \Nat(F_1 \otimes F_2, G_1 \otimes G_2)
  \end{equation*}
  induced from the tensor product is injective.
\end{lemma}
\begin{proof}
  Let $f_1, \cdots, f_m \in \Nat(F_2, G_2)$ be linearly independent elements. Suppose that $c_1, \cdots, c_m \in \Nat(F_1, G_1)$ are elements such that
  \begin{equation*}
    c_1|_X \otimes f_1|_Y + \cdots + c_m|_X \otimes f_m|_Y = 0
  \end{equation*}
  for all $(X, Y) \in \mathcal{I}_1 \times \mathcal{I}_2$. If we fix $X \in \mathcal{I}_2$, we can apply Lemma~\ref{lem:faithful-2} and obtain that $c_i|_X = 0$ for $i = 1, \cdots, m$. By letting $X$ run through all objects of $\mathcal{I}_1$, we have that $c_i = 0$ for $i = 1, \cdots, m$. Thus the map under consideration is injective.
\end{proof}

Now we can prove Theorem~\ref{thm:lin-indep}. Our proof is based on a proof of the linear independence of group\-like elements in a coalgebra over a field.

\begin{proof}[Proof of Theorem~\ref{thm:lin-indep}]
  Recall the assumptions: $\mathcal{I}$ is a skeletally small monoidal category, $\mathcal{C}$ is a tensor category over $\bfk$ and $F$ and $G$ are strong monoidal functors from $\mathcal{I}$ to $\mathcal{C}$.

  We first note that $0 \not\in \Nat_\otimes(F, G)$. Indeed, if $g \in \Nat_\otimes(F, G)$, $g|_\mathbf{1}: F(\mathbf{1}) \to G(\mathbf{1})$ must be an isomorphism. Since $F(\mathbf{1}) \cong \mathbf{1} \ne 0$, $g|_\mathbf{1} \ne 0$, and hence $g \ne 0$.

  Suppose to the contrary that $\Nat_\otimes(F, G)$ is linearly dependent. Then there exist elements $g_1, \cdots, g_m \in \Nat_\otimes(F, G)$ and $\lambda_1, \cdots, \lambda_m \in \bfk$ such that
  \begin{equation*}
    g := \lambda_1 g_1 + \cdots + \lambda_m g_m \in \Nat_\otimes(F, G)
  \end{equation*}
  and $g \ne g_i$ for $i = 1, \cdots, m$. We may suppose that $g_1, \cdots, g_m$ are linearly independent. Since $g \ne 0$, we may also suppose that $\lambda_h \ne 0$ for some $h$.

  By the definition of monoidal natural transformations, we have
  \begin{equation*}
    \sum_{j = 1}^m \lambda_j g_j|_X \otimes g_j|_Y = g|_X \otimes g|_Y
    = \sum_{i, j = 1}^m \lambda_i \lambda_j g_i|_X \otimes g_j|_Y.
  \end{equation*}
  for all $X, Y \in \mathcal{I}$. By Lemma~\ref{lem:faithful-3}, $\sum_{i = 1}^m \lambda_j \lambda_i g_i = \lambda_j g_j$ for each $j$. By the linear independence of $g_i$'s, we have that $\lambda_i = 0$ for $i \ne h$. Therefore, $g = g_h$. This is a contradiction.
\end{proof}

\subsection{Proof of Corollary~\ref{cor:fin-nat}}

In view of Theorem~\ref{thm:lin-indep}, it is interesting to know when $\Nat(F, G)$ is fi\-nite-di\-men\-sion\-al. We conclude this section by proving the following lemma, which completes the proof of Corollary~\ref{cor:fin-nat}.

\begin{lemma}
  \label{lem:nat-fin-dim}
  Let $F, G: \mathcal{A} \to \mathcal{B}$ be right exact $\bfk$-linear  functors between $\bfk$-linear abelian categories. Suppose that $\mathcal{A}$ has a projective generator $P$ and every object of $\mathcal{A}$ is of finite length. Then the linear map
  \begin{equation*}
    (-)|_P: \Nat(F, G) \to \Hom_\mathcal{B}(F(P), G(P)), \quad f \mapsto f|_P
  \end{equation*}
  is injective. In particular, if every hom-set in $\mathcal{B}$ is fi\-nite-di\-men\-sion\-al, then
  \begin{equation*}
    \dim_\bfk \Nat(F, G) \le \dim_\bfk \Hom_\mathcal{B}(F(P), G(P)) < \infty.
  \end{equation*}
\end{lemma}

Note that a finite $\bfk$-linear abelian category satisfies the assumptions on $\mathcal{A}$ in this lemma. Therefore Corollary~\ref{cor:fin-nat} follows immediately from Theorem~\ref{thm:lin-indep}.

\begin{proof}
  Let $f \in \Nat(F, G)$. By the assumption, for every $X \in \mathcal{A}$, there exists an exact sequence $P^{\oplus m} \to P^{\oplus n} \to X \to 0$. By applying $F$ and $G$ to this sequence, we have a commutative diagram
  \begin{equation*}
    \begin{CD}
      F(P)^{\oplus m} @>>> F(P)^{\oplus n} @>>> F(X) @>>> 0 \\
      @V{(f|_P)^{\oplus m}}VV @V{(f|_P)^{\oplus n}}VV @V{f|_X}VV \\
      G(P)^{\oplus m} @>>> G(P)^{\oplus n} @>>> G(X) @>>> 0
    \end{CD}
  \end{equation*}
  in $\mathcal{B}$ with exact rows. This means that $f|_X$ is determined by $f|_P$, and hence the map $f \mapsto f|_P$ is injective.
\end{proof}

\section{Some remarks on $\Aut_\otimes(\id_\mathcal{C})$}
\label{sec:structure-auto}

\subsection{Bound of the order}

Let $\mathcal{C}$ be a finite tensor category over $\bfk$. In this section, we give some remarks on the structure of the group $G(\mathcal{C}) = \Aut_\otimes(\id_\mathcal{C})$. We first note that by the definition of natural transformations, $G(\mathcal{C})$ is abelian.

Let $I$ be the set of isomorphism classes of simple objects of $\mathcal{C}$. For each $i \in I$, we fix $S_i \in i$. Let $P_i$ be the projective cover of $S_i$. Then $P = \bigoplus_{i \in I} P_i$ is a projective generator. Applying Lemma~\ref{lem:nat-fin-dim}, we have an injective linear map $\Nat(\id_\mathcal{C}, \id_\mathcal{C}) \to \End_\mathcal{C}(P)$. By the definition of natural transformations, the image of this map must be contained in the center of $\End_\mathcal{C}(P)$. Hence we obtain an injective linear map
\begin{equation*}
  \Nat(\id_\mathcal{C}, \id_\mathcal{C}) \to Z(\End_\mathcal{C}(P)), \quad \eta \mapsto \eta|_{P}.
\end{equation*}

\begin{remark}
  As every object of $\mathcal{C}$ is of finite-length, $\mathcal{C}$ is $\bfk$-linearly equivalent to the category of fi\-nite-di\-men\-sion\-al right $\End_\mathcal{C}(P)$-modules, and hence the above map is actually an isomorphism of $\bfk$-algebras. We omit the detail here since the surjectivity will not be used.  
\end{remark}

Theorem~\ref{thm:lin-indep} states that $G(\mathcal{C}) \subset \Nat(\id_\mathcal{C}, \id_\mathcal{C})$ is linearly independent. Therefore, we have the following bound on the order of $G(\mathcal{C})$.

\begin{proposition}
  $|G(\mathcal{C})| \le \dim_\bfk Z(\End_\mathcal{C}(P))$.
\end{proposition}

In the case where $\mathcal{C} = \Rep(H)$ for some fi\-nite-di\-men\-sion\-al Hopf algebra $H$, this proposition is obvious since the right-hand side of the inequality is equal to the dimension of the center of $H$.

\subsection{Values on simple objects}

Let us consider the map
\begin{equation*}
  \Nat(\id_\mathcal{C}, \id_\mathcal{C}) \to \prod_{i \in I} \End_\mathcal{C}(S_i),
  \quad \eta \mapsto (\eta|_{S_i})_{i \in I}.
\end{equation*}
This map is not injective any more unless $\mathcal{C}$ is semisimple. Since $\bfk$ is algebraically closed, we can identify $\End_\mathcal{C}(S_i)$ with $\bfk$. The above map induces a group homomorphism
\begin{equation*}
  \varphi: G(\mathcal{C}) \to \prod_{i \in I} \End_\mathcal{C}(S_i)^\times \cong \Map(I, \bfk^\times),
  \quad g \mapsto (i \mapsto g|_{S_i}).
\end{equation*}

We show that $\varphi$ is injective if $\bfk$ is of characteristic zero. To describe the kernel of $\varphi$, we introduce some subgroups of $G(\mathcal{C})$. If $p = \mathop{\rm char}(\bfk) > 0$, then we set
\begin{align*}
  G(\mathcal{C}) _p & = \{ g \in G(\mathcal{C}) \mid \text{$g^{p^k} = 1$ for some $k \ge 0$} \}, \\
  G(\mathcal{C})'_p & = \{ g \in G(\mathcal{C}) \mid \text{The order of $g$ is relatively prime to $p$} \}.
\end{align*}
Otherwise we set $G(\mathcal{C})_p = \{ 1 \}$ and $G(\mathcal{C})'_p = G(\mathcal{C})$. By the fundamental theorem of finite abelian groups, we have a decomposition $G(\mathcal{C}) = G(\mathcal{C})_p \times G(\mathcal{C})_p'$.

\begin{lemma}
  \label{lem:auto-indec-1}
  Let $X \in \mathcal{C}$ be an indecomposable object.
  \begin{enumalph}
  \item If $g \in G(\mathcal{C}) _p$, then $g|_X$ is unipotent.
  \item If $g \in G(\mathcal{C})'_p$, then $g|_X = \lambda \cdot \id_X$ for some $\lambda \in \bfk^\times$.
  \end{enumalph}
\end{lemma}
\begin{proof}
  Let $g \in G(\mathcal{C})$. As $X$ is indecomposable, $\End_\mathcal{C}(X)$ is a local algebra. Since $\bfk$ is algebraically closed, $g|_X$ can be written uniquely in the form
  \begin{equation}
    \label{eq:def-lambda}
    g|_X = \lambda \id_X + r \quad (\lambda \in \bfk^\times, r \in \mathfrak{m}),
  \end{equation}
  where $\mathfrak{m}$ is the maximal ideal of $\End_\mathcal{C}(X)$. If $r \ne 0$, then there exists $k \ge 1$ such that $r \in \mathfrak{m}^{k-1}$ but $r \ne \mathfrak{m}^{k}$. By the binomial formula, we have
  \begin{equation*}
    (g|_X)^n \equiv \lambda^n \id_X + n r \pmod{\mathfrak{m}^{k}}
  \end{equation*}
  for every $n \ge 0$. This implies that if $(g|_X)^n = \id_X$, then $\lambda^n = 1$ and $n r = 0$.

  {\rm (a)} Suppose that $g \in G(\mathcal{C})_p$. Since the claim is obvious for $p = 0$, we assume that $p > 0$. Then the order of $g|_X$ is a power of $p$. Thus, by the above observation, $\lambda = 1$. This implies that $g|_X$ is unipotent.

  {\rm (b)} Suppose that $g \in G(\mathcal{C})'_p$. Then the order of $g|_X$ is nonzero in $\bfk$. Thus, by the above observation, $r = 0$. This implies that $g|_X = \lambda \id_X$.
\end{proof}

In what follows, we denote $\lambda$ in equation~(\ref{eq:def-lambda}) by $\lambda_g(X)$. We remark the following easy but important property of $\lambda_g(X)$.

\begin{lemma}
  \label{lem:auto-indec-2}
  Let $X$ and $Y$ be indecomposable objects of $\mathcal{C}$. If $X$ and $Y$ belong to the same block, then $\lambda_g(X) = \lambda_g(Y)$ for every $g \in G(\mathcal{C})'_p$.
\end{lemma}

Here, a {\em block} is an equivalence class of indecomposable objects of $\mathcal{C}$ under the weakest equivalence relation such that two indecomposable objects $X$ and $Y$ of $\mathcal{C}$ are equivalent whenever $\Hom_\mathcal{C}(X, Y) \ne 0$.

\begin{proof}
  Let $g \in G(\mathcal{C})_p'$. By the definition of blocks, it is sufficient to prove in the case when there exists a nonzero morphism $f: X \to Y$. By the naturality of $g$,
  \begin{equation*}
    \lambda_g(X) f = (g|_X) \circ f = f \circ (g|_Y) = \lambda_g(Y) f
  \end{equation*}
  Since $f \ne 0$, we have $\lambda_g(X) = \lambda_g(Y)$.
\end{proof}

Now we have a description of the kernel of $\varphi$ as follows:

\begin{proposition}
  \label{prop:auto-kernel}
  $\Ker(\varphi) = G(\mathcal{C})_p$. In particular, $\varphi$ is injective if $p = 0$.
\end{proposition}
\begin{proof}
  Let $g \in G(\mathcal{C})_p$. Then, by Lemma~\ref{lem:auto-indec-1}, $g|_S = \id_S$ for every simple object $S$ of $\mathcal{C}$ and hence $g \in \Ker(\varphi)$. This implies that $G(\mathcal{C})_p \subset \Ker(\varphi)$.

  Next let $g \in G(\mathcal{C})'_p \cap \Ker(\varphi)$. Let $X$ be an indecomposable object of $\mathcal{C}$ and fix a simple sub\-object $S$ of $X$. Then, by Lemma~\ref{lem:auto-indec-2}, $\lambda_g(X) = \lambda_g(S)$. On the other hand, $\lambda_g(S) = 1$ since $g \in \Ker(\varphi)$. This implies that $g|_X = \id_X$ for all indecomposable object $X$ and hence $g = 1$. Therefore $G(\mathcal{C})'_p \cap \Ker(\varphi) = \{1\}$.

  Now we recall that $G(\mathcal{C}) = G(\mathcal{C})_p \times G(\mathcal{C})'_p$. Our claim follows immediately from the above observations.
\end{proof}

It is interesting to characterize the image of $\varphi$. Let $N_{i j}^k$ ($i, j, k \in I$) be the multiplicity of $S_k$ as a composition factor of $S_i \otimes S_j$. In the case where $\mathcal{C}$ is semisimple, it is known that the image of $\varphi$ is the set of functions $\lambda: I \to \bfk^\times$ such that
\begin{equation}
  \label{eq:image-1}
  \text{$\lambda(i) \lambda(j) = \lambda(k)$ whenever $N_{i j}^k \ne 0$}
  \quad (i, j, k \in I).
\end{equation}
In general, the image of $\varphi$ is smaller than the set of such functions.

\begin{proposition}
  \label{prop:auto-image}
  The image of $\varphi$ is the set of all functions satisfying~(\ref{eq:image-1}) and the following condition:
  \begin{equation}
    \label{eq:image-2}
    \text{$\lambda(i) = \lambda(j)$ whenever $S_i$ and $S_j$ belong to the same block} \quad (i, j \in I).
  \end{equation}
\end{proposition}
\begin{proof}
  We remark that if an indecomposable object $X \in \mathcal{C}$ has $S_i$ as a composition factor, then $S_i$ and $X$ belong to the same block. Indeed, then there exists a nonzero morphism $P_i \to X$ and hence $P_i$ and $X$ belong to the same block. On the other hand, since $S_i$ is a quotient of $P_i$, $S_i$ and $P_i$ belong to the same block. Therefore the claim follows.

  Let $\lambda = \varphi(g)$. By Proposition~\ref{prop:auto-kernel}, we may assume $g \in G(\mathcal{C})_p'$. (\ref{eq:image-2}) follows from Lemma~\ref{lem:auto-indec-2}. Thus we check that $\lambda$ satisfies~(\ref{eq:image-1}). Let $i, j \in I$. By the definition of monoidal natural transformations,
  \begin{equation*}
    g|_{S_i \otimes S_j} = g|_{S_i} \otimes g|_{S_j} = \lambda(i) \lambda(j) \id_{S_i \otimes S_j}.
  \end{equation*}
  Suppose that $N_{i j}^k \ne 0$. This means that $S_i \otimes S_j$ has $S_k$ as a composition factor. Let $X$ be an indecomposable direct summand of $S_i \otimes S_j$ having $S_k$ as a composition factor. By the above equation, $g|_X = \lambda(i) \lambda(j) \id_X$. On the other hand, since $X$ and $S_k$ belong to the same block, $g|_X = \lambda(k)\id_X$. Therefore $\lambda(k) = \lambda(i) \lambda(j)$.

  Conversely, given a function $\lambda: I \to \bfk^\times$ satisfying~(\ref{eq:image-1}) and~(\ref{eq:image-2}), we define a natural automorphism $g: \id_\mathcal{C} \to \id_\mathcal{C}$ as follows: If $X$ is an indecomposable object of $\mathcal{C}$, then $g|_X = \lambda(i)\id_X$ where $i \in I$ is such that $X$ has $S_i$ as a composition factor. As $\lambda$ satisfies~(\ref{eq:image-2}), this does not depends on the choice of $i$. We can extend $g$ to all objects of $\mathcal{C}$, since they are direct sums of indecomposable objects.

  Now we need to show that $g \in G(\mathcal{C})$, that is, $g|_{X \otimes Y} = g|_X \otimes g|_Y$ for all objects $X, Y \in \mathcal{C}$. We may assume that $X$ and $Y$ are indecomposable. Suppose that
  \begin{equation*}
    X = \sum_{i \in I} m_i S_i \quad (m_i \in \mathbb{Z}_{\ge 0}) \text{\quad and \quad}
    Y = \sum_{j \in I} n_j S_j \quad (n_i \in \mathbb{Z}_{\ge 0})
  \end{equation*}
  in the Grothendieck ring $K(\mathcal{C})$. Then
  \begin{equation*}
    X \cdot Y = \sum_{k \in I} ( \sum_{i, j \in I} m_i n_j N_{i j}^k ) S_k
  \end{equation*}
  in $K(\mathcal{C})$. This equation means that if $X \otimes Y$ has $S_k$ as a composition factor, then there exist $i, j \in I$ such that $X$ has $S_i$ as a composition factor, $Y$ has $S_j$ as a composition factor and $N_{i j}^k \ne 0$. By the definition of $g$ and~(\ref{eq:image-1}), we have $g|_{X \otimes Y} = g|_X \otimes g|_Y$.

  It is obvious that $\lambda = \varphi(g)$. The proof is completed.
\end{proof}

The following theorem is a direct consequence of Proposition~\ref{prop:auto-kernel} and~\ref{prop:auto-image}.

\begin{theorem}
  If $p = 0$, then $\varphi$ gives an isomorphism between $G(\mathcal{C})$ and the group of functions $\lambda: I \to \bfk^\times$ satisfying~(\ref{eq:image-1}) and~(\ref{eq:image-2}).
\end{theorem}

\section*{Acknowledgements}

The author is supported by Grant-in-Aid for JSPS Fellows.


\end{document}